\documentclass{amsart}
\usepackage{mathptmx, amssymb, colonequals}

\newtheorem{theorem}{Theorem}
\theoremstyle{definition}
\newtheorem{conjecture}[theorem]{Conjecture}
\newtheorem{example}[theorem]{Example}

\def\ge{\geqslant}

\def\ann{\operatorname{ann}}
\def\image{\operatorname{image}}

\begin{document}
\title{On a conjecture of Lynch}

\author{Anurag K. Singh}
\address{Department of Mathematics, University of Utah, 155 South 1400 East, Salt Lake City, UT~84112, USA}
\email{singh@math.utah.edu}

\author{Uli Walther}
\address{Department of Mathematics, Purdue University, 150 N University St., West Lafayette, IN~47907, USA}
\email{walther@math.purdue.edu}

\thanks{A.K.S.was supported~by NSF grant DMS~1801285, and U.W.~by the Simons Foundation Collaboration Grant for Mathematicians~\#580839}
\maketitle


The following conjecture has recently attracted attention, e.g.,~\cite{Boix:Eghbali1, Boix:Eghbali2, DSZ, Hochster:Jeffries}:

\begin{conjecture}~\cite[Conjecture~1.2]{Lynch}
Let $R$ be a local ring, and $I$ an ideal of $R$. If the cohomological dimension of $I$ is a positive integer $c$, then
\[
\dim R/\ann_R H^c_I(R) = \dim R/H^0_I(R).
\]
\end{conjecture}

The conjecture is know to be false: the first counterexamples were constructed in~\cite{Bahmanpour}; these are nonequidimensional, with $\dim R\ge 5$. We present here a modification---with a short, elementary proof---that serves as a counterexample with $\dim R=3$. This is a counterexample, as well, to \cite[Proposition~4.3]{Lynch} and to \cite[Theorem~4.4]{Lynch}; the error there is in the chain of inequalities in the proof of Proposition~4.3, \cite[page~550]{Lynch}, in the reduction from a complete local ring to a complete local unmixed ring: the cohomological dimension may change under the reduction step.

\begin{example}
Let $k$ be a field, and set $R\colonequals k[x,y,z_1,z_2]/(xyz_1,\ xyz_2)$. Consider the local cohomology module $H^2_{(x,y)}(R)$. Using a \v Cech complex on $x$ and $y$, one sees that
\[
H^2_{(x,y)}(R) = R_{xy}/\image(R_x+R_y).
\]
The images of $z_1$ and $z_2$ are zero in $R_{xy}$, so the local cohomology module above agrees with~$S_{xy}/\image(S_x+S_y)$, where $S\colonequals R/(z_1, z_2)$ is isomorphic to the polynomial ring~$k[x,y]$. Hence $H^2_{(x,y)}(R)$ is a nonzero $R$-module, with annihilator $(z_1, z_2)$. On the other hand, since the ideal $(x,y)$ contains nonzerodivisors, $H^0_{(x,y)}(R)=0$. Hence one has
\[
\dim R/\ann_R H^2_{(x,y)}(R) = 2\qquad\text{whereas}\qquad \dim R/H^0_{(x,y)}(R)=3.
\]
\end{example}


\end{document}